\documentclass{article}
\usepackage{amsmath,amsfonts,amsthm,graphics}
\graphicspath{{figures/}}

\newtheorem{theorem}{Theorem}
\newtheorem{lemma}[theorem]{Lemma}

\newtheorem{corollary}[theorem]{Corollary}

\newcommand\bbp{{\widetilde{\mathfrak X}}}
\newcommand\bp{{\mathfrak X}}
\newcommand\rw{{\widetilde{\mathfrak X}}}
\newcommand\rwX{{\widetilde X}}
\newcommand\bpt{\rw^2}
\newcommand\bpo{\rw^1}
\newcommand\bphi{{\widetilde\phi}}
\newcommand\Z{{\mathbb Z}}
\newcommand\DD{{\mathcal D}}
\newcommand\PP{{\mathcal P}}
\newcommand\RR{{\mathbb R}}
\newcommand\TT{{\mathbb T}}
\newcommand\NN{{\mathbb N}}
\newcommand\eps{\varepsilon}

\newcommand\dd{\,d}
\newcommand\norm[1]{\ensuremath{\|#1\|}}
\newcommand\bb[1]{\bigl(#1\bigr)}
\newcommand\area{\mathrm{area}}

\newcommand{\refT}[1]{Theorem~\ref{#1}}
\newcommand{\refC}[1]{Corollary~\ref{#1}}
\newcommand{\refL}[1]{Lemma~\ref{#1}}

\newcommand\E{{\mathop{\mathbb E{}}\nolimits}}
\renewcommand\Pr{{\mathop{\mathbb P{}}\nolimits}}
\newcommand\Po{\operatorname{Po}}

\newcommand\Gnop{G_{n,\omega,p}}
\newcommand\Gope{G_{\omega,p,\eps}}
\newcommand\Gnol{G_{n,\omega,\la/\omega}}
\newcommand\Gnolp{G_{n,\omega,\la'/\omega}}
\newcommand\bGop{{\widetilde G}_{\omega,p}}
\newcommand\Gop{G_{\omega,p}}
\newcommand\Gdrop{G_{d,r,\omega,p}}
\newcommand\Gol{G_{\omega,\la/\omega}}
\newcommand\Gle{G_{\la,\eps}}
\newcommand\Glee{G_{\la/\eps,\eps}}
\newcommand\GGlee{G(\la/\eps,C_\eps)}
\newcommand\pc{{p_{\mathrm c}}}
\newcommand\lac{{\la_{\mathrm c}}}
\newcommand\ac{{a_{\mathrm c}}}
\newcommand\la{\lambda}
\newcommand\Ga{\Gamma}
\newcommand\rhoh{{\rho_{\mathrm h}}}
\newcommand\rhov{{\rho_{\mathrm v}}}

\begin{document}
\title{Line-of-sight percolation}
\date{January 28, 2007; revised April 7, 2008}

\author{B\'ela Bollob\'as\thanks{Department of Mathematical Sciences,
University of Memphis, Memphis TN 38152, USA}
\thanks{Trinity College, Cambridge CB2 1TQ, UK}
\thanks{Research supported in part by NSF grants CCR-0225610,
 DMS-0505550 and W911NF-06-1-0076}
\and Svante Janson%
\thanks{Department of Mathematics, Uppsala University,
 PO Box 480, SE-751 06 Uppsala, Sweden}
\and Oliver Riordan%
\thanks{Royal Society Research Fellow, Department of Pure Mathematics
and Mathematical Statistics, University of Cambridge, Wilberforce Road,
Cambridge CB3 0WB, UK}}
\maketitle

\begin{abstract}
Given $\omega\ge 1$, let $\Z^2_{(\omega)}$ be the graph with vertex set $\Z^2$ in
which two vertices are joined if they agree in one coordinate and
differ by at most $\omega$ in the other. (Thus $\Z^2_{(1)}$ is precisely
$\Z^2$.) Let $\pc(\omega)$ be the critical probability for site percolation
on $\Z^2_{(\omega)}$. Extending recent results of Frieze, Kleinberg, Ravi and Debany,
we show that $\lim_{\omega\to\infty} \omega\pc(\omega)=\log(3/2)$. We
also prove analogues of this result for the $n$-by-$n$ grid and in higher
dimensions, the latter involving interesting connections to Gilbert's
continuum percolation model.
To prove our results, we explore the component of the origin in a certain
non-standard way, and show that this exploration is well approximated
by a certain branching random walk.
\end{abstract}

\section{Introduction and results}

Frieze, Kleinberg, Ravi and Debany~\cite{FKRD}
proposed the following random graph as a model of an {\em
ad hoc} network in an environment with (regular) obstructions.  Given
positive integers $n$ and $\omega$, let $G=[n]^2_{(\omega)}$ be the graph with vertex
set $[n]\times [n]$ in which two vertices $(x_1,y_1)$ and $(x_2,y_2)$
are joined if $x_1=x_2$ and $|y_1-y_2|\le \omega$ or $y_1=y_2$ and
$|x_1-x_2|\le \omega$.  Let $V$ be a random subset of $[n]\times [n]$
obtained by selecting each point $(x,y)$ with probability $p$,
independently of the other points, and let $G[V]$ denote the subgraph
of $G$ induced by the vertices in $V$.
We write $\Gnop$ for $G[V]$, and $G_{\omega,p}$ for the infinite
random subgraph of $\Z^2_{(\omega)}$
defined in the same way but starting from $\Z^2$
rather than $[n]\times [n]$.

The model $\Gnop'$ studied
in~\cite{FKRD} was defined on the $n$ by $n$ torus rather than
the grid, but $\Gnop$ seems more natural in terms of the original motivation.
For the results here, it will make no difference which variant we consider.

One interpretation of $\Gnop$ is as follows: sensors are dropped onto
a random selection of crossroads in a regularly laid out city (or
planted forest); two sensors can
communicate if they are within distance $\omega$ and the line of sight
between them is not blocked by a building (tree). The graph $\Gnop$
indicates which pairs of sensors can communicate directly, so we would
like to know when $\Gnop$ has a giant component, and roughly how large
it is. Alternatively, taking the point of view of percolation
theory, we would like
to study percolation in $\Gop$, i.e., to know for which choices
of the parameters $\Gop$ has an infinite component.

Kolmogorov's 0/1--law implies that for each fixed $\omega$
there is a `critical probability' $\pc(\omega)$ such that if
$p>\pc(\omega)$ then $\Gop$ contains an infinite component
with probability $1$,
while if $p< \pc(\omega)$ then with probability $1$
all components of $\Gop$ are finite.
Furthermore, by the `uniqueness theorem' of Aizenman, Kesten and Newman~\cite{AKNunique}
(see also~\cite[p.~121]{BRbook}), when an infinite cluster exists it is unique
with probability 1.
Among other results, Frieze, Kleinberg, Ravi and Debany~\cite{FKRD}
proved (essentially) that
\[
 1/(4e) \le \liminf_{\omega\to\infty}\omega\pc(\omega)
 \le \limsup_{\omega\to\infty}\omega\pc(\omega) \le 1,
\]
and posed the natural question of determining the value
of $\lim \omega\pc(\omega)$, assuming it exists. (In fact,
they worked with finite random graphs, as we shall
in \refT{th_2} below.)
Here we shall answer this question, proving that
\begin{equation}\label{lim}
 \lim_{\omega\to\infty} \omega\pc(\omega) =\log(3/2)=0.4054\dots,
\end{equation}
as part of the more detailed \refT{th_1} below.

In homogeneous cases such as this, where the expected degrees
of all vertices are equal, the critical probability
is perhaps more naturally described in terms of the `critical
expected degree'. Since the expected degree of a vertex of
$\Gop$ is $4p\omega$, \eqref{lim} states that the critical
expected degree of $\Gop$
tends to $4\log(3/2)=1.6218\dots$ as $\omega\to\infty$.

Note that it is too much to hope to find
$\pc(\omega)$ exactly for a given $\omega$: with $\omega$ fixed,
$\pc(\omega)$ is the critical probability for site percolation on a certain
lattice $\Z^2_{(\omega)}$ (non-planar for $\omega\ge 2$); the exact values 
of such critical probabilities are known only in a few
very special cases. Indeed, $\pc(1)$ is the critical probability
for site percolation on $\Z^2$, and even this is not known
exactly.

In spirit, \refT{th_1} and its proof are very close to the arguments
given in~\cite{BJRspread} for
Penrose's spread-out percolation model~\cite{Penrose}:
when $\omega\to\infty$ with $p\omega$ bounded, we can
approximate the neighbourhood
exploration process in $\Gop$ by a certain branching process,
and describe the distribution of the small components in this graph
in terms of the branching process. We then show by a rather {\em ad hoc}
argument that almost all vertices not in `small' components
are in a single infinite component. An analogous assertion holds
for $\Gnop$ and its unique giant component.
In contrast to \cite{BJRspread}, where we could appeal to the general
sparse inhomogeneous random graph model of~\cite{BJRkernels}, here
even the local coupling to the branching process has to be done by hand. We shall
discuss this later.

To state our results precisely, we first define the branching process
alluded to above. Let $\mu>0$ be a real parameter, and let $\Ga=\Ga_\mu$
denote (a random variable with) the geometric distribution with parameter
$1-e^{-\mu}$. Thus,
\[
 \Pr(\Ga=k) = (1-e^{-\mu})^k e^{-\mu}
\]
for $k=0,1,\ldots$, and $\E(\Ga)=e^{\mu}-1$.
Let $\Ga_\mu^{(r)}$ denote the distribution of the sum of $r$ independent
copies of $\Ga_\mu$.

Let $\bp_\mu$ be the branching process $(X_0,X_1,\ldots)$
in which $X_0$ consists of a single particle $x_0$,
each particle $x\in X_t$ has children independently of the other particles
and of the history, the number of children of $x_0$ has the distribution
$\Ga_\mu^{(4)}$,
and the number of children of a particle $x\in X_t$, $t\ge 1$,
has the distribution $\Ga_\mu^{(2)}$.
Let $\phi(\mu)$ denote the {\em survival probability} of the branching process $\bp_\mu$,
i.e.,
\[
 \phi(\mu) = \Pr\bb{\hbox{$|X_t|>0$ for all $t$}}.
\]
Note that
\[
 \E\bb{\Ga_\mu^{(2)}} = 2\E(\Ga_\mu) = 2(e^{\mu}-1),
\]
so the branching process $\bp_\mu$
is supercritical if and only if $\mu>\log(3/2)$, 
i.e., $\phi(\mu)>0$ if and only if $\mu>\log(3/2)$.

In the result below, $C_0$ denotes the component of the origin
in the random graph $\Gop$; if the origin is not a vertex
of this graph, then we set $C_0=\emptyset$.
We write $\Pr_{\omega,p}$ for the probability measure
associated to $\Gop$, and
$\theta(\omega,p)=\Pr_{\omega,p}(|C_0|=\infty)$
for the probability that $C_0$ is infinite.

\begin{theorem}\label{th_1}
Let $0<\la<\log(3/2)$ be constant. If $\omega$ is large enough
then all components of $\Gol$ are finite with probability $1$.
Furthermore, there is a constant $a_\la>0$ such that,
for all large enough $\omega$,
\[
 \Pr_{\omega,\la/\omega}(|C_0|\ge k) \le \frac\la\omega e^{-a_\la (k-1)}
\]
holds for every $k\ge 1$.

Let $\la>\log(3/2)$ be constant.
If $\omega$ is large enough, then $\Gol$
contains a unique infinite component with probability $1$. Furthermore,
\begin{equation}\label{thform}
 \theta(\omega,\la/\omega) \sim \phi(\la)\la/\omega
\end{equation}
as $\omega\to\infty$.
\end{theorem}

Note that Theorem~\ref{th_1} certainly implies \eqref{lim}, which is equivalent
(as noted earlier) to the statement that the critical expected degree tends
to $4\log(3/2)$ as $\omega\to\infty$.
The proof of Theorem~\ref{th_1} will form the bulk of Section~\ref{sec_proofs}.
In fact, the upper bound on component sizes (giving the lower bound on the
critical expected degree), is fairly easy. The lower bound on component
sizes (giving the upper bound on the critical expected degree) will be
derived in a way that will enable us to prove a finite analogue
of Theorem~\ref{th_1}, stated as Theorem~\ref{th_2} below.

For the upper bound on component sizes
we shall prove a much more detailed statement than stated in Theorem~\ref{th_1}, namely that
\begin{equation}\label{exactupper}
 \Pr_{\omega,p}(|C_0|\ge k) \le p\Pr(|\bp_{\mu}|\ge k)
\end{equation}
for every $\omega\ge 1$, $0<p<1$ and $k\ge 1$, where
\begin{equation}\label{laform}
 e^{-\mu} = (1-p)^\omega,
\end{equation}
and $|\bp_\mu|$ denotes the total number of particles in all generations
of $\bp_\mu$. 
For $\mu<\log(3/2)$, the branching process $\bp_\mu$ is strictly subcritical,
and standard results imply that there is some $b_\mu>0$ such that
$\Pr(|\bp_\mu|\ge k)\le e^{-b_\mu (k-1)}$ holds for every $k$.
If $\la<\log(3/2)$ is constant
and we set $p=\la/\omega$, then as $\omega\to\infty$ we have
$\mu\sim p\omega=\la$. In particular, $\mu$ is bounded above
by some $\mu'<\log(3/2)$ when $\omega$ is large enough, so the first statement
in \refT{th_1} follows.

Letting $k$ tend to infinity in \eqref{exactupper} (or simply setting
$k=\infty$, which makes perfect sense), we see that
$\theta(\omega,p)\le p\phi(\mu)$. Standard results on branching processes
imply that $\phi(\cdot)$ is continuous. 
Thus \eqref{exactupper} implies the upper bound implicit
in \eqref{thform}, i.e., that
\[
 \theta(\omega,\la/\omega) \le (1+o(1))\phi(\la)\la/\omega
\]
as $\omega\to\infty$ with $\la>\log(3/2)$ constant.

Note that for the exact comparison with a branching process above, it is natural
to work with the parameter $\mu=-\log((1-p)^\omega)$.
However, in the lower bounds on component
sizes, we can obtain only approximate results, and it is more natural
to use the asymptotically equivalent parameter $\lambda=p\omega$, which
we have used throughout the statement of \refT{th_1}. We shall switch
between these two parameters as and when convenient.

Although our main focus is the infinite random graph $\Gop$,
our other aim is to prove a result corresponding to Theorem~\ref{th_1} for the finite
graphs $\Gnop$ and $\Gnop'$.
Note that there is a natural coupling in which $\Gnop$, defined on the grid,
is a subgraph of $\Gnop'$, defined on the torus.
Let $C_1(G)$ denote the number of vertices
in a largest component of a graph $G$.

\begin{theorem}\label{th_2}
Let $0<\la<\log(3/2)$ be constant. There are constants $A_\la$
and $\omega_\la$ such that
\[
 C_1(\Gnol)\le C_1(\Gnol')\le A_\la\log n
\]
holds with probability $1-o(1)$ as $n\to\infty$,
whenever $\omega=\omega(n)$ satisfies $\omega_\la\le\omega\le n$.

Let $\la>\log(3/2)$ be constant. If $\omega=\omega(n)$ is such
that $\omega$ and $n/\omega\to\infty$, then
\[
 C_1(\Gnol),\  C_1(\Gnol') = (\phi(\la)+o_p(1))n \la/\omega
\]
as $n\to\infty$.
\end{theorem}

We shall also state results for the case $\omega\to\infty$ with $n=\Theta(\omega)$;
for details see Subsection~\ref{noc}.

\section{Proofs}\label{sec_proofs}

\subsection{Upper bounds on component sizes}
We start with some simple upper bounds on the component sizes.
Let $p$ and $\omega$ be arbitrary. To prove \eqref{exactupper},
it suffices to couple $\Gop$ and $\bp_\mu$
so that if we explore the component $C_0$ of the origin
in $\Gop$ in a suitable manner,
then the number of vertices reached at each
stage is at most the number of particles of $\bp_\mu$
in the corresponding generation,
where $\mu=-\omega\log(1-p)$ is given by \eqref{laform}.

In each step of our overall exploration process (except for the first)
we either `explore vertically', first up and then down,
or `explore horizontally', first to the left and then
to the right.
To `explore upwards' from a vertex $v=(x,y)$ of $\Gop$,
we test the points $(x,y+1)$, $(x,y+2)$, \ldots, one by one. Each
test `succeeds' if the relevant point is a vertex of $\Gop$ that was
not previously reached (in an earlier stage of the overall exploration).
We stop the upwards exploration as soon as
$\omega$ consecutive tests fail.

If $v_1,\ldots,v_r$
denote the vertices of $\Gop$ reached during this upwards exploration,
corresponding to the successful tests,
then $vv_1$, $v_1v_2,\ldots,v_{r-1}v_r$ are all edges of $\Gop$.
Also, there is typically no edge of $\Gop$ from $v_r=(x,y')$ to a point
$(x,y'')$ with $y''>y'$; this is not always true as a test may fail
due to the presence of a vertex we reached earlier.
If we have not previously tested the relevant
points, so each is present in $\Gop$ with probability $p$, 
independently of the others, then the probability that the next
$\omega$ tests fail is exactly
$(1-p)^\omega = e^{-\mu}$, so $r$ has the geometric distribution
$\Ga_\mu$. In general, we may have previously tested some of the points;
a repeated test of a point always fails by definition, so, conditional
on previous steps of the exploration, the distribution of $r$
is stochastically dominated by $\Ga_\mu$.

To `explore vertically' from a vertex $v$, we explore upwards from
$v$, and then explore downwards from $v$. Similarly, to explore
horizontally, we first explore to the left, and then to the right.

Let us condition on $v_0=(0,0)$ being a vertex of $\Gop$, an event
of probability $p$. Set $Y_0=\{v_0\}$. Explore vertically and horizontally
from $v_0$, and let $Y_1$ be the set of vertices of $\Gop$
reached in these explorations. From the remarks above, $|Y_1|$
has exactly the distribution $\Ga_\mu^{(4)}$.
Suppose that we have defined $Y_t$, $t\ge 0$.
For every vertex $v\in Y_t$
that was reached during a horizontal exploration, we explore vertically
from $v$. Similarly, we explore horizontally from each $v\in Y_t$ reached
during a vertical exploration. Let $Y_{t+1}$ be the set of new vertices
of $\Gop$ reached.

It is not hard to check that the exploration just defined does indeed uncover
the whole component $C_0$. Indeed, let $C_0'=\bigcup_{t\ge 0} Y_t$.
Then we certainly have $v_0\in C_0'\subset C_0$. Since $C_0$ is a connected
subgraph of $\Gop$, if $C_0'\ne C_0$ then there are vertices
$v\in C_0'$ and $w\in C_0\setminus C_0'$ that are adjacent in $\Gop$.
Choosing $v$ and $w$ at minimal (Euclidean) distance, there are no vertices of $\Gop$
between $v$ and $w$. We may suppose without loss of generality
that $v=(x,y)$ and $w=(x',y)$ with $0<x'-x\le \omega$. When exploring, 
we reach $v$ at some stage; if we do so during a vertical exploration,
then exploring from $v$ to the right we reach $w$ at the first step.
On the other hand, if we reach $v$ while exploring to the right or to
the left, the same exploration would also reach $w$, either just after
or just before reaching $v$. This shows that $w\in C_0'$, contradicting
our assumptions.

From the remarks above, conditional on previous explorations, the number of vertices
reached in a horizontal or vertical exploration is stochastically dominated
by the $\Ga_\mu^{(2)}$ distribution. Thus the whole sequence $|Y_0|,|Y_1|,|Y_2|,\ldots$
is stochastically dominated by $|X_0|,|X_1|,|X_2|,\ldots,$ where $X_0,X_1,\ldots$ are the generations
of the branching process $\bp_\mu$.
Since $C_0=\bigcup_{t\ge 0} Y_t$, we thus have
\[ 
 \Pr\bb{|C_0|\ge k \mid v_0\in V(\Gop)} \le \Pr(|\bp_\mu|\ge k)
\]
for every $k\ge 1$, which gives \eqref{exactupper}.

As noted earlier, the first part of \refT{th_1} follows, along with the 
upper bound implicit in \eqref{thform}.

\subsection{The lower bound on component sizes}

Throughout this section we shall set $p=\la/\omega$, where $\la>\log(3/2)$ is
constant and we let $\omega\to \infty$. All asymptotic notation refers
to this limit.

Usually, in showing that the neighbourhood exploration process in some
random graph is well approximated by a branching process, one first
shows that a given vertex is unlikely to be in a short cycle.
Here the latter statement is false: given that $v=(x,y)$ is a vertex of $\Gop$,
it has probability $\Theta(1)$ of being in a triangle of the form
$(x,y)$, $(x,y+a)$, $(x,y+b)$, $0<a<b\le\omega$.
However, these triangles (and similar cycles) have already been accounted
for in the construction of the branching process.
If we delete all edges of $\Gop$ corresponding
to line segments whose interiors contain vertices of $\Gop$,
then it is easy to check that the density of short cycles in 
the resulting graph is very low.
We shall show something more or less equivalent to this statement.

Let $k$ be fixed. As before, let us condition on the event that
$v_0=(0,0)$ is a vertex
of $\Gop$. We claim that we can couple the first $k$ steps
$Y_0,Y_1,\ldots,Y_k$ of the exploration defined above with the first $k$
generations of $\bp_\la$ so that they agree (in the sense that
$|Y_t|=|X_t|$ for $0\le t\le k$)
with probability $1-o(1)$. (Here the asymptotic notation refers to the limit
$\omega\to\infty$ with $k$ and $\la=p\omega$ fixed.)
We have already shown that the exploration is stochastically dominated
by $\bp_\mu$, where $\mu=-\omega\log(1-p)$. Since $\mu\sim\la$,
the processes $\bp_\mu$ and $\bp_\la$ may be coupled so that their first
$k$ generations agree whp, so it suffices to `bound $(Y_t)$ from below',
i.e., to couple a subset of this exploration process with $\bp_\la$
so as to agree whp for $k$ generations.

\begin{figure}[ht]
 \centering
 \input{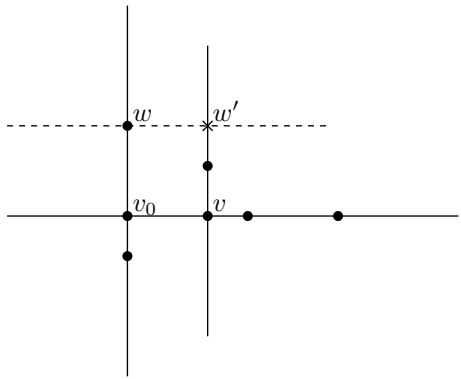}
 \caption{Part of our exploration, including a possible clash (between
$w$ and $w'$) that we must avoid. The solid lines indicate
regions that have already been tested for vertices of $\Gop$;
these lines extend a distance $\omega$ from the last vertex
found.}\label{fig_1}
\end{figure}
In doing so, we must be careful: suppose that as part of our exploration
we find vertices of $C_0$ located as in Figure~\ref{fig_1}. Then we are in trouble --
since $w$ and $w'$ lie on the same horizontal line, the horizontal explorations
from these points will necessarily interfere. In particular,
the expected number of new vertices reached from $w$ and $w'$ 
together will be significantly less than
$2\E(\Ga_\la^{(2)})$.
Thus, if we reach $w$ first, we do not allow ourselves to test
the point $w'$.

In general, we consider the {\em lines} of $\Z^2$, i.e., the sets
of the form $\{x\}\times \Z$ or $\Z\times \{y\}$.
Initially, we mark all lines as `unclaimed',
except the two lines through $v_0$,
which we mark as `claimed by $v_0$'.
Whenever we find a vertex $v\in C_0$ during
a horizontal exploration, this vertex immediately claims the vertical
line it lies on.
Similarly, when we find a vertex $v$ during a vertical exploration, $v$
claims the horizontal line it lies on. We modify our exploration
as follows: when exploring horizontally or vertically from
a vertex $v$, we omit testing any point $w$ on a line
already claimed by a vertex other than $v$. This ensures
that the same line cannot be claimed twice: when we
find a vertex $w\in C_0$ by exploring from $v$, exactly one line through
$w$ has been claimed, namely that joining $v$ to $w$; the vertex
$w$ then claims the other unclaimed line.

In this modified exploration, a point $(x,y)$ can only be tested once,
so each tested point is present in $\Gop$ with (conditional)
probability $p$. Since the expected number of vertices reached within
$k$ steps of the exploration is $O(1)$, whp there are at most
$\log\omega$ such vertices, say, and hence at most $2\log\omega$
claimed lines.
In each horizontal or vertical exploration, we are exploring
along a claimed line, and we omit testing
at most one point on each other claimed line.
Thus 
we omit testing $O(\log\omega)=o(1/p)$ points. Since we would have expected
only $o(1)$ successes among
these tests, omitting them makes essentially no difference: the number $r$
of new vertices found may be coupled with a $\Ga_{\la}^{(2)}$ distribution
so as to agree with probability $1-o(1)$. As we expect
to reach $O(1)$ vertices during 
the first $k$ steps of
our exploration, the sum of these $o(1)$
error probabilities is still $o(1)$, so the claim follows.

The truth of the statement `$|C_0|\ge k$' is certainly determined by
the numbers $|Y_0|,\ldots,|Y_{k-1}|$.
Thus, for every $1\le k<\infty$ we have
\[
 \Pr_{\omega,\la/\omega}\bb{|C_0|\ge k\mid v_0\in V(\Gol)} \to
 \Pr(|\bp_\la|\ge k)
\]
as $\omega\to\infty$, i.e.,
\[
 \Pr_{\omega,\la/\omega}(|C_0|\ge k) \sim  \Pr(|\bp_\la|\ge k)\la/\omega.
\]
It follows that there is some $K=K(\omega)\to\infty$ such that
\[
 \Pr_{\omega,\la/\omega}( |C_0|\ge K) \sim \Pr(|\bp_\la|\ge K)\la/\omega.
\]
Since the branching process $\bp_\la$ does not depend on $\omega$, the
right-hand side above is asymptotically $\phi(\la)\la/\omega$. 
To complete the proof of \eqref{thform}
it remains to show that
\[
 \Pr_{\omega,\la/\omega}(K\le|C_0|<\infty) = o(\la/\omega).
\]
In other words, roughly speaking, we must show that almost all vertices of
$\Gol$ in `large' components are in infinite components.
As noted in the introduction, the probability
that there is an infinite component is either 0 or 1, and when
there is an infinite component it is unique with probability 1,
so the remaining assertion of \refT{th_1} follows.

\bigskip
So far we have coupled the neighbourhood exploration process
in $\Gop$ with a branching process $\bp_\la$. If we keep track
of the locations of the points in $Y_0,Y_1,\ldots$ as well as their
number, then the appropriate limit object is a branching random walk.

Indeed, let us turn $\bp_\la$ into a branching random walk $\rw_\la$
on $\RR^2$
(or, more formally, on $\RR^2\times\{\mathrm{v},\mathrm{h}\}$) as follows. Each
generation $\rwX_t$, $t\ge 1$, will consist of a finite set of points of
$\RR^2$, each labelled with either `v' (for reached by a vertical step)
or `h' (reached by a horizontal step).  A particle of type `v' at a
point $(x,y)$ has children in the next generation according to the
following rule: generate a Poisson process $\PP$ of intensity $\la$ on
$\RR\times\{y\}$. Starting from $(x,y)$ and working to the right,
include as children of $(x,y)$ all points of $\PP$ until the first time
that we come to a gap of length greater than $1$. Do the same to the
left. All children have type `h'. The rule for children of particles
of type `h' is similar, using the line $\{x\}\times \RR$ instead.
Sometimes we start with two particles at the origin, one of type `h' and one
of type `v'; we write $\bpt_\la$ for this branching random walk.
At other times, we consider the same branching rule but start
with a single particle, writing $\bpo_\la$. When the starting
rule is clear or unimportant, we write simply $\rw_\la$.
Since the number of offspring of a particle in $\rw_\la$ is independent
of the position of this particle, we may view $\bpt_\la$
as our original branching process $\bp_\la$ by simply ignoring the positions
of the particles.

It is easy to check that the rules above correspond to a certain
limit of the horizontal and vertical explorations we considered earlier,
where we take $\omega\to\infty$ and rescale by dividing the coordinates
of our lattice points by $\omega$.
In particular, the coupling argument above shows that for any fixed $k$,
we may couple $(Y_t)_{t\le k}$ with the first $k$ generations $(\rwX_t)_{t\le k}$
of the
branching random walk $\bpt_\la$ so that with probability $1-o(1)$ there
is exactly one particle $(x',y')\in \rwX_t$
for each $(x,y)\in Y_t$, $1\le t\le k$,
and $|x/\omega-x'|,|y/\omega-y'|\le k/\omega$.

We shall only study the branching random walk $\rw_\la$ in a trivial way,
calculating certain expectations; moreover, we shall avoid all
detailed calculations, needing only `soft' arguments.

We start with a simple lemma showing that the supercritical branching
random walk $\rw_\mu$, $\mu>\log(3/2)$, remains supercritical
even when restricted to a large enough square. By the {\em restriction}
of $\rw_\mu$ to a region $R$ we mean the branching random walk
obtained from $\rw_\mu$ by deleting all particles
that lie outside $R$, along with all their descendants.
We write $\rw_\mu(C)$ for the restriction of $\rw_\mu$ to
$[-C,C]^2$.

\begin{lemma}\label{l_sq}
Let $\mu>\log(3/2)$ and $\eps>0$ be fixed. There are constants $A$, $B$ and $T_1$
(depending on $\mu$)
with the following properties.

(i) For any point $(x_0,y_0)$ of $[-A,A]^2$, if we start the 
branching random walk $\rw_\mu$ with a single point (of type `h', say)
at $(x_0,y_0)$ and restrict to $[-B,B]^2$, then the expected
number of points of generation $T_1$ inside $[-A,A]^2$ is at least $2$.

(ii) We have
\begin{equation}\label{surv}
 \Pr\bb{ \bpt_\mu(B) \hbox{ survives forever}} \ge \phi(\mu)-3\eps,
\end{equation}
where $\phi(\mu)$ is the survival probability of the branching process $\bp_\mu$.
\end{lemma}

\begin{proof}
Starting the unrestricted walk $\bpo_\mu$ with a single
particle at the origin, the expected
size of generation $t$ is $\E\bb{\Ga_\mu^{(2)}}^t=(2e^{\mu}-2)^t$.
Since $\mu>\log(3/2)$, there is some $T_1$ such that this expectation
is at least $10$. Fix such a $T_1$.

Let $\bpo_\mu(A)$ denote the restriction of $\bpo_\mu$ to $[-A,A]^2$.
As $A\to\infty$, the
first $T_1$ generations of $\bpo_\mu(A)$ converge in distribution
to the first $T_1$ generations of $\bpo_\mu$. Hence there is some $A_0$ such
that the expected size of generation $T_1$ of $\bpo_\mu(A)$
is at least $8$ whenever $A\ge A_0$.

By symmetry, the distributions of the number of points of generation
$T_1$ of $\bpo_\mu(A_0)$ in each of the four squares $[0,\pm A_0]\times [0,\pm A_0]$
are identical. Let $\DD$ denote this distribution, noting that $\E(\DD)\ge 2$.

Let $A\ge A_0$ be a constant depending on $\mu$ and $\eps$, to be chosen later.
Set $B=2A$, and let $(x_0,y_0)$ be any point of $[-A,A]^2$.
Without loss of generality, we may assume that $x_0,y_0\le 0$, 
so $[x_0,x_0+A]\times [y_0,y_0+A]$ is contained in $[-A,A]^2$;
see Figure~\ref{fig_2}.
\begin{figure}[ht]
 \centering
 \input{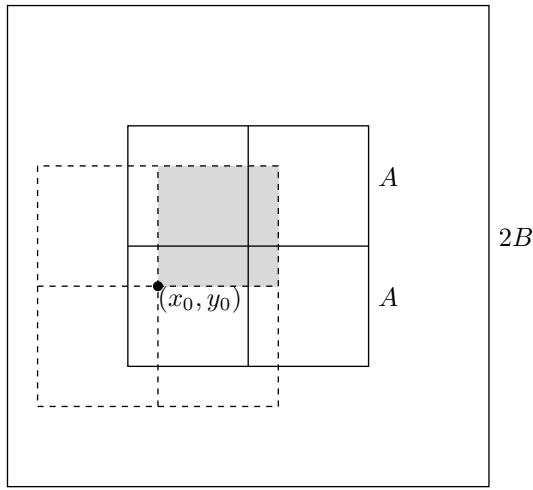}
 \caption{The dotted lines show the square $S$ of side $2A$ centered
at $(x_0,y_0)$, which is contained within the outer
square of side $2B=4A$ centered at the origin. The shaded quadrant
of $S$ is contained in $[-A,A]^2$.}
\label{fig_2}
\end{figure}
Let $\rw_\mu'$ denote the branching random walk $\bpo_\mu$
started at $(x_0,y_0)$
and restricted to $[-B,B]^2$.
Since the square of side $2A$ centered at $(x_0,y_0)$ is contained in
$[-B,B]^2$, the process $\rw_\mu'$ stochastically dominates
the translated process $\bpo_\mu(A)+(x_0,y_0)$.
In particular, the number of points of $\rw_\mu'$ in $[-A,A]^2$,
which is at least the number of points in $[x_0,x_0+A]\times [y_0,y_0+A]$,
stochastically dominates $\DD$. Since $\E(\DD)\ge 2$, this proves (i).

We have shown that, for any $A\ge A_0$,
if we start our branching walk restricted
to $[-2A,2A]^2$ at any point of $[-A,A]^2$, and look at the number of points
of $[-A,A]^2$ reached at times $0$, $T_1$, $2T_1$, $3T_1$,\ldots,
then the numbers we see stochastically dominate a certain supercritical
Galton--Watson branching process (with offspring distribution
$\DD$). This process survives forever with some
probability $p_0>0$ depending on $\mu$ only.

Recall that ignoring the positions of the particles turns
the (unrestricted) branching walk $\bpt_\mu$ into our original
branching process $\bp_\mu$. Except that the offspring distribution
is different for the first generation, $\bp_\mu$ is a standard Galton--Watson
process, so standard results (for example, Theorem 2 in Athreya and Ney~\cite{AN})
imply that with probability $1$, either $\bp_\mu$ dies out
(an event of probability $1-\phi(\mu)$), or $|X_t|\to\infty$ as $t\to\infty$.
It follows that there is a time $T$
such that with probability
at least $\phi(\mu)-\eps$ generation $T$ contains at least $\log(1/\eps)/p_0$
particles. 
Hence, choosing $A\ge A_0$ large enough, the walk $\bpt_\mu(A)$
has probability at least $\phi(\mu)-2\eps$ of generating at least
$\log(1/\eps)/p_0$ particles in generation $T$; these particles
lie in $[-A,A]^2$ by definition.

Let $N$ denote the number of particles of $\bpt_\mu(B)$
in generation $T$ in $[-A,A]^2$. Since $\bpt_\mu(B)$
may be regarded as a superset of $\bpt_\mu(A)$,
we have $N\ge \log(1/\eps)/p_0$ with probability
at least $\phi(\mu)-2\eps$.
Each particle of $\bpt_\mu(B)$ in $[-A,A]^2$ survives forever with probability
at least $p_0$, independently of the others, so
\begin{eqnarray*}
 \Pr\bb{\bpt_\mu(B)\hbox{ dies out}}
  &\le& \Pr\bb{N<\log(1/\eps)/p_0} + (1-p_0)^{\log(1/\eps)/p_0} \\
 &<& 1-(\phi(\mu)-2\eps)+\eps = 1-(\phi(\mu)-3\eps).
\end{eqnarray*}
In other words, $\bpt_\mu(B)$ survives with
probability at least $\phi(\mu)-3\eps$, proving \eqref{surv}.
\end{proof}

In the next lemma we consider the restricted walk $\rw_\mu(B)$
started at a particle $v$ of type `h' just below the bottom
edge of $[-B,B]^2$; making sense of this requires a slight modification
of our definition of restriction.
Recalling that we label particles by the kind of exploration that reached
them, the first step from $v$ is vertical, so the first generation
$\rwX_1$
of the unrestricted walk started at $v$ consists of points on the same
vertical line as $v$. To define the first generation $\rwX_1'$
of our restriction $\rw_\mu(B)$, first delete all points of $\rwX_1$
outside $[-B,B]^2$. Then, if some point 
in the remaining set $\rwX_1\cap[-B,B]^2$
is within distance
$1$ of $v$, let $\rwX'_1=\rwX_1\cap[-B,B]^2$;
otherwise,
delete all remaining points too and let $\rwX'_1=\emptyset$.
All points in $\rwX'_1$ are
inside
$[-B,B]^2$, and from here we continue the restricted walk as usual.
The reason for the somewhat fussy definition of the first generation
is that when we return to the graph $\Gop$, we will wish to find certain
paths starting at a vertex $v$ just outside $[-\omega B,\omega B]^2$,
with all remaining vertices inside $[-\omega B,\omega B]^2$;
the first edge of such a path must join $v$ to a point in $[-\omega B,\omega B]^2$
at distance at most $\omega$ from $v$.

\begin{figure}[ht]
  \centering
 \input{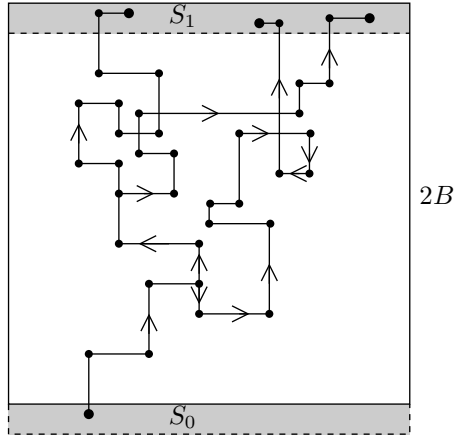}
 \caption{The shaded regions are the sets $S_0$ and $S_1$ considered
in Lemma~\ref{l_B} and Corollary~\ref{c_cont}; the solid square
is the boundary of $[-B,B]^2$. The paths are for $T_2=16$.
Note that there may be additional points of the branching process
in the interiors of the line segments of these paths; indeed,
there must be in any line segment of length more than $1$.} \label{f_3}
\end{figure}
\begin{lemma}\label{l_B}
Let $\mu>\log(3/2)$ be given. There are constants $B$ and $T_2$
(with $B$ as in \refL{l_sq})
with the following property. Whenever we start $\rw_\mu(B)$ at a point
of type `h' in the set $S_0=[-B,B]\times [-B-1/2,-B]$ shown
in Figure~\ref{f_3},
the expected number of points of generation $T_2$ of type `h' in the set
$S_1=[-B,B]\times [B-1/2,B]$ 
is at least $2$.
\end{lemma}
\begin{proof}
Let $A$, $B$ and $T_1$ be as in Lemma~\ref{l_sq},
and let $(\rwX_t')$ denote the generations of
the branching random walk under consideration.
From the definition of $\rw_\mu$, the (non-normalized)
probability density function
describing the $x$- or $y$- displacement in a single step
is positive on the whole real line.
Hence the probability that $\rwX_2'$
contains a point of $[-A,A]^2$, which depends on $(x_0,y_0)$,
is bounded below by some constant $c=c(B)>0$, independent of
$(x_0,y_0)$. 
Using part (i) of \refL{l_sq} and induction on $k$, it follows that
\[
 \E\bb{ |\rwX_{kT_1+2}' \cap [-A,A]^2| } \ge 2^k c
\]
for every $k$.

Arguing as at the start of the proof,
there is some $c'>0$ such that for any point $(x,y)\in [-A,A]^2$,
if we start a copy of $\bpo_\mu$ at a point $(x,y)$ of type `h'
and restrict to $[-B,B]^2$,
the expected number of points of $S_1$ that we obtain two generations later is
at least $c'$.
It follows that
\[
 \E\bb{ |\rwX_{kT_1+4}' \cap S_1| } \ge 2^k cc'.
\]
Choosing $k$ large enough that $2^kcc'>2$, and taking $k$
even so that all points of generation $T_2=kT_1+4$ have
type `h', the result follows.
\end{proof}

The same proof shows that \refL{l_B} still holds if, instead of
looking for points of type `h' in $S_1$, we look for points of
type `v' in the vertical strip $S_1'=[B-1/2,B]\times [-B,B]$ 
just inside the right hand edge of $[-B,B]^2$.
The same is true in the following two corollaries.

\begin{corollary}\label{c_cont}
Let $\mu>\log(3/2)$ and $\eps>0$ be given, and let
$B$ and $T_2$ be as in \refL{l_B}. Then every large enough constant $N$
has the following property. Whenever we start $\rw_\mu$
at $N$ points of $S_0$ of type `h' and restrict to $[-B,B]^2$,
the probability that generation $T_2$ contains at least $N$ points
of $S_1$ with distinct $x$-coordinates is at least $1-\eps$.
\end{corollary}
\begin{proof}
Let $B$ and $T_2$ be as in Lemma~\ref{l_B}, and
let $Z_i$ be the number of $T_2$th generation descendants of the $i$th
starting point that lie in $S_1$. Since we start $N$ independent
(restricted) branching walks, the random variables $Z_i$ are independent;
\refL{l_B} shows that each has expectation at least $2$.
Although the $Z_i$ are not identically distributed,
as $T_2$ is constant we have a common upper bound (of order
$(2e^{\mu}-2)^{2T_2}$) on $\E(Z_i^2)$ for any $i$.
Hence $Z=\sum_{i=1}^N Z_i$ has expectation at least $2N$
and variance $O(N)$.
Consequently, $\Pr(Z\ge N)\ge 1-\eps$
if $N$ is large enough.

Recalling that $T_2$ is even, in going from generation $T_2-1$ to $T_2$
we take horizontal steps, so the $x$-coordinates of the points we reach
are distinct with probability $1$.
\end{proof}

The supercriticality of the restricted walk and the
survival probability bound \eqref{surv} have
the following consequence.

\begin{corollary}\label{c_start}
Let $\mu>\log(3/2)$ and $\eps>0$ be fixed, let $B$ be as above, and
let $N$ be a constant. There is a constant $T_3$ such that, 
if $\bpt_\mu(B)$ is started (as usual) with two particles at $(0,0)$,
one of type `h' and one of type `v', then with probability at least
$\phi(\mu)-4\eps$ either generation $T_3$ or generation $T_3+1$
contains at least $N$ points of $S_1$
with type `h' with distinct $x$-coordinates.
\end{corollary}
\begin{proof}
This follows easily from Lemma~\ref{l_sq}(ii),
using arguments
similar to those in the proofs of \refL{l_B} and Corollary~\ref{c_cont}. 
We omit the details.
Note that we need to consider two generations, $T_3$ and $T_3+1$, as it may be that
only one of the initial particles survives, in which case, after a while,
either only odd generations or only even generations contain particles of type `h', and
we cannot say in advance which it will be.
\end{proof}

We are now ready to complete the proof of \refT{th_1}.
\begin{proof}[Proof of \refT{th_1}]
We have already proved the upper bounds on component sizes; it remains
to prove the lower bound implicit in \eqref{thform}, i.e., to show that
\begin{equation}\label{lb}
 \liminf_{\omega\to\infty} \theta(\omega,\la/\omega)(\la/\omega)^{-1} \ge \phi(\la)
\end{equation}
whenever $\la>\log(3/2)$.
As remarked above, it follows by standard results
that if $\omega$ is large enough (so $\theta(\omega,\la/\omega)>0$),
then with probability 1 there is a unique
infinite component.

Let $\la>\log(3/2)$ be fixed, and set $\mu=\la$ for compatibility
with our branching process notation above.
(Often, at this point in the proof
one would fix an arbitrary $\mu<\la$, but here the usual `elbow room'
turns out not to be needed.)

Fix $\eps>0$ with $10\eps<\phi(\mu)$,
and let $B$, $N$, $T_2$ and $T_3$ be defined as in the corollaries above.
As before, we condition on $(0,0)$ being a vertex of $\Gop$,
and let $(Y_t)$ denote the neighbourhood exploration process in $\Gop$,
$p=\mu/\omega$.
We have already seen that, for any fixed $k$,
the first $k$ generations of the rescaled process $(Y_t/\omega)$ 
(defined by $Y_t/\omega=\{(x/\omega,y/\omega):(x,y)\in Y_t\}$)
may be coupled with the branching random walk $\bpt_\mu=(\rwX_t)$
so that with probability $1-o(1)$ as $\omega\to\infty$
the first $k$ generations agree up to displacements
of individual points by up to $k/\omega$.

This coupling result holds also if we explore only the subgraph
of $\Gop$ induced by vertices in $[-\omega B,\omega B]^2$,
and replace $\bpt_\mu$ by its restriction $\bpt_\mu(B)$ to $[-B,B]^2$.
Combining this observation with \refC{c_start} above,
we see that if $\omega$ is large enough then, with probability
at least $\phi(\mu)-5\eps$, exploring $\Gop$ from $(0,0)$
within $[-\omega B,\omega B]^2$, in either $T_3$ or $T_3+1$
steps we reach at least $N$ vertices
$v_1,\ldots,v_N$
of $\omega S_1=[-\omega B,\omega B]\times [ \omega(B-1/2),\omega B]$,
with the last step to each vertex being horizontal.
Note that we may assume that $v_1,\ldots,v_N$ have distinct $x$-coordinates
(since the corresponding particles of $\rw_\mu$ do).

For $(a,b)\in \Z^2$, let $Q_{a,b}$ be the square
\[
 Q_{a,b} = \bb{(2a-1)\omega B,(2a+1)\omega B}
  \times \bb{(2b-1)\omega B,(2b+1)\omega B}.
\]
Since we exclude the boundaries, the squares $Q_{a,b}$ are disjoint.
Since $N$ is constant (i.e., independent of $\omega$), the coupling
above extends to explorations started at any $N$ points 
of $\Z^2$,
provided that there are no clashes at the first step, i.e.,
that we do not attempt to explore vertically from two points
on the same vertical line, or horizontally from two points on the
same horizontal line. In what follows, we can always
assume that this proviso is satisfied (from the `distinct $x$-coordinates'
conclusions of Corollaries~\ref{c_cont} and~\ref{c_start}, and
the corresponding `distinct $y$-coordinates' conclusions of the variants
where we end with points near the right of $[-B,B]^2$). We shall
not comment on this annoying technicality further. 

Suppose that, without testing points in $Q_{a,b}$, we have found
$N$ vertices $v_1,\ldots,v_N$ of $\Gop$
lying just below $Q_{a,b}$, i.e., in the strip
of width $\omega/2$ bordering $Q_{a,b}$ from below.
Then, from \refC{c_cont}, provided $\omega$ is large enough,
with probability at least $1-2\eps$ we may find paths in $\Gop$
from $\{v_i\}$ to $N$ points $w_i$ of $Q_{a,b}$ just below the top
of $Q_{a,b}$, using only vertices in $Q_{a,b}$.
Similarly (using the variant of \refC{c_cont} with $S'_1$ in place
of $S_1$), with probability
at least $1-2\eps$ we may find such 
paths to $N$ points just inside the right-hand side of $Q_{a,b}$.
By symmetry, the same conclusions hold if
we  start with $N$ points $v_i$ just to the left of $Q_{a,b}$.

\begin{figure}[htb]
 \centering
 \input{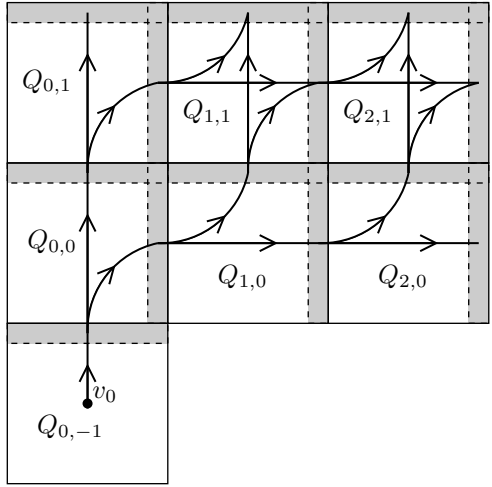}
 \caption{A schematic depiction of the proof that percolation occurs. 
Roughly speaking, whenever
there is an arrow joining two shaded regions, from $N$ vertices in the first
shaded region we are very likely to find $N$ in the second. We start by
showing that from $v_0$ we find $N$ vertices in the shaded region
of $Q_{0,-1}$ with probability close to $\phi(\mu)$.}\label{f_perc}
\end{figure}

To complete the proof, it is convenient to explore the component $C_0$
of $\Gop$ containing the point $v_0=(0,-2B\omega)\in Q_{0,-1}$,
in a manner shown schematically in Figure~\ref{f_perc}.
Recalling that $v_0$ is present in $\Gop$ with probability $p$,
by \refC{c_start},
with unconditional probability at least $p(\phi(\mu)-5\eps)$
we find paths within $\Gop\cap Q_{0,-1}$ from $v_0$
to $N$ points just below the top of $Q_{0,-1}$. Let us denote
this
event by $E$.

Let us now explore within the squares $Q_{a,b}$, $a,b\ge 0$,
working upwards in layers, i.e., in increasing order of $a+b$.
When we come to a particular square $Q_{a,b}$, if we have found
$N$ vertices of $\Gop$ near the top of the square $Q_{a,b-1}$ below,
or $N$ vertices near the right of the square $Q_{a-1,b}$,
then with probability at least $1-4\eps$, testing only points
inside $Q_{a,b}$, we find paths from these vertices to $N$
points near the right of $Q_{a,b}$, and to $N$ points near the top
of $Q_{a,b}$. If this happens, we declare the point $(a,b)\in \NN^2$
to be {\em open}. If we did not reach $N$ vertices of $\Gop$ near the top
of the square below or the right of the square to the left,
we declare $(a,b)$ to be open regardless of what happens inside $Q_{a,b}$.

Because each exploration is confined to its own square,
each $(a,b)$ is declared open with conditional probability at least $1-4\eps$.
Thus the distribution
of open points $(a,b)\in \NN^2$ stochastically dominates a distribution
in which each point is open independently with probability $1-4\eps$.

Suppose that $E$ holds and that there is an infinite
oriented path $P$ in $\NN^2$ starting at $(0,0)$, in which
every vertex is open and every edge goes up or to the right.
Then $v_0$ lies
in an infinite component of $\Gop$. 
It is easy to check that the probability of the existence of such a path $P$
tends to $1$ as $\eps\to 0$.
(For example, one can bound the expected number
of `blocking cycles' as in the proof of Lemma 6 in~\cite{BRbook},
and show that this expectation tends to $0$ as $\eps\to 0$.)
Thus there is a function $f(\eps)$
with $f(\eps)\to 0$ as $\eps\to 0$ such that
\[
 \theta(\omega,\mu/\omega) \ge (1-f(\eps))\Pr(E) \ge (1-f(\eps)(\mu/\omega)(\phi(\mu)-5\eps)
\]
whenever $\omega$ is large enough.
Since $\eps>0$ was arbitrary,
this proves \eqref{lb}.
\end{proof}

We now turn to the finite graphs $\Gnop$ and $\Gnop'$.
\begin{proof}[Proof of \refT{th_2}]
Let $v=(x,y)$, $1\le x,y\le n$, be any potential
vertex of $\Gnop$, and let $C_v$ and $C_v'$ denote
the components of $v$ in the graphs $\Gnop$ and $\Gop$.
We may regard $\Gnop$ as a subgraph of $\Gop$, so
\[
 \Pr(|C_v|\ge k) \le \Pr(|C_v'|\ge k).
\]
If $\la<\log(3/2)$ is constant and $\omega$ is large enough, then
by \refT{th_1} the second
probability above is at most $(\la/\omega)e^{-a_\la (k-1)}$ for every $k\ge 1$,
where $a_\la>0$ depends on $\la$ only.
Setting $A=3/a_\la$, say, taking $k=A\log n$, and recalling that there
are $n^2$ choices for $v$, it follows that
with probability $1-o(1)$ the graph $\Gnol$ contains no components
of size larger than $A\log n$, whenever $\omega$ is large enough.

Turning to $\Gnop'$, although this cannot be seen as a subgraph of $\Gop$,
for any given vertex $v$ of $\Gnop'$ we may couple the neighbourhood exploration
processes in $\Gnop'$ and in $\Gop$ started at $v$ so that the latter dominates.
(Alternatively, note that the branching process domination argument we gave
in $\Gop$ applies just as well in $\Gnop'$.)
Thus the bound above holds for $\Gnop'$ also. Since $\Gnop$ is a subgraph
of $\Gnop'$, the first statement in \refT{th_2} follows.

\medskip
It remains to show that in the supercritical case, if $\la>\log(3/2)$
is constant and $\omega=\omega(n)$ with $\omega\to\infty$ and $n/\omega\to\infty$,
then $\Gnol$ and $\Gnol'$ contain a unique giant component of the `expected'
size, namely $(\phi(\la)+o(1))n\la/\omega$.
As the arguments are mostly rather standard, we shall only sketch them.

The local coupling arguments above show that we have the expected number
of vertices in `small' components: it remains only to show that whp there
is a component of size at least $(\phi(\la)+o(1))n\la/\omega$.
Fix $\eps>0$ and set $\la'=\la-\eps$; we assume that $\eps$ is such that
$\la-2\eps>\log(3/2)$. 
It is also convenient to assume that $\la<1$.
We start by analyzing components
in $\Gnolp$, using the extra vertices of $\Gnol$ for later sprinkling.

We know from \refL{l_sq} that if $C$ is a large enough constant (depending
on $\la'$ and $\eps$), then the branching random walk $\rw_{\la'-\eps}(C)$ restricted
to $[-C,C]^2$ is supercritical, and indeed survives with probability
at least $\phi(\la'-\eps)-\eps$ whenever we start it with two particles of
type `h' and `v' at a point not too
close to the boundary of $[-C,C]^2$. Let $G$ be the subgraph of $\Gnolp$
induced by vertices in $[0,2\omega C]^2$. 
When we explore the component of a vertex of $G$ not too close to the boundary,
there are two sources of error in coupling this exploration to the corresponding
(restricted) branching random walk: the first comes from approximating the
discrete distribution of points on a horizontal or vertical line by
a Poisson process, which gives rise to an error probability of $O(\la'/\omega)=O(1/\omega)$
at each step. The second comes from omitting tests for points on lines
already claimed by other points; we may account
for this by reducing the branching process
parameter slightly. If we have found $m$ points so far,
each claims only one line (apart from the first). It follows
that for some $L=\Theta(\eps\omega)$,
we can couple our exploration to dominate $\rw_{\la'-\eps}(C)$ as
long as $m\le L$, with an error
probability of at most $\la'/\omega$ in each step
(coming from the approximation of a binomial distribution by a Poisson), 
and thus a total error probability that is $O(\eps)$.

Let $N$ be the number
of vertices of $G$ that are in components of size at least $L$.
From the remarks above it follows that
$\E N \ge xN_0$, where
$x=\phi(\la'-\eps)-O(\eps)=\phi(\la)-O(\eps)$, 
and $N_0=(2C\omega)^2\la'/\omega$ is the expected
number of vertices of $G$.
Furthermore, the number $N$ is concentrated in the sense that
$\E|N-\E N|=O(\eps N_0)$ (at least for large $\omega$);
this can be seen, for example, 
by comparing with $N_k$, the number of vertices of $G$ in components
of order at least $k$, for some fixed but large $k$, and noting that
starting our exploration from two different vertices and using
a corresponding upper bound, we find that the variance of $N_k$ is
$o(N_0^2)$. 
Hence, $N>(x-\eps^{1/2})N_0$ with probability at least $1-O(\eps^{1/2})$.

Sprinkling extra vertices with density $(\la-\la')/\omega$, it is easy to check that, whp,
all components in $G$ of order at least $L$
join up to form a single component,
which thus has size at least $N$. One crude argument is as follows:
suppose that before sprinkling we have components $C_1$ and $C_2$ each with
at least $L=\Theta(\omega)$ vertices. Let $A$ denote the set
of $x$-coordinates of vertices of $C_1$, and $B$ the set
of $y$-coordinates of vertices of $C_2$.
Since the expected number of vertices of $G$
an any horizontal or vertical line in $[0,2\omega C]^2$ is $O(1)$,
whp no line contains more than $O(\log \omega)$ vertices, so we may assume
that $|A|,|B| =\Omega(\omega/\log \omega)$, say.
Let $S$ be the set of sprinkled vertices in $A\times B$. Then $|S|$
is concentrated about its mean, which is of order at least $\omega/(\log \omega)^2$,
so we may assume that $|S|\ge \omega/(\log \omega)^3$, say. Using
once more the fact that there aren't too many vertices on any line, we may find
a subset $S'$ in which all coordinates are distinct, with $|S'|=\Omega(\omega/(\log \omega)^5)$,
say. Finally, for each $(x,y)\in S'$, we look for a vertical path
of sprinkled vertices (outside $A\times B$) joining $(x,y)$ to the (nearest)
vertex of $C_1$ with the same $x$-coordinate, and a horizontal path to $C_2$.
We find such paths if each of at most $8C$ intervals of length $\omega/2$
contains a sprinkled vertex, an event whose probability is bounded away
from zero. These events are independent for different $(x,y)\in S'$,
so with very high probability at least one such pair of paths is present.

Returning to $\Gnol$, or $\Gnol'$,
we may cover the vertex set of this graph by squares of side-length $2C\omega$, overlapping
in regions of width $C\omega$, say. The argument above shows that within each square,
we almost certainly find a component of the right size, i.e., containing
at least a proportion $x-\eps^{1/2}$ of the vertices. The sprinkling above also shows
that the giant components of overlapping squares are very likely to meet.
Considering either mixed or $1$-independent percolation on a $K$ by $K$ grid
in which each vertex/bond is open with probability at least $p$, by estimating
the expected number of blocking cycles it is easy to check that the probability
that any two vertices are joined is at least $f(p)$, for some function $f(p)$
tending to $1$ as $p\to 1$, independent of $K$. Comparison with such a model
shows that as $n,n/\omega\to\infty$, whp the large components
in almost all squares are linked up, giving a giant component in
$\Gnol$ or $\Gnol'$
containing at least a proportion $\phi(\la)-\delta(\eps)$ of all
vertices, for some function $\delta(\eps)\to0$ as $\eps\to0$.
\end{proof}

Note that the condition $n/\omega\to\infty$ in the second part of
\refT{th_2} is essential. Indeed, if $n/\omega$ is bounded, 
say $n\le C\omega$, then the approximation
of the component exploration by the branching random walk $\rw_\la$ breaks
down even during the first step:  this branching random walk has probability
$\Theta(1)$ of starting by moving to a point at distance at least $C$ (corresponding
to distance at least $C\omega\ge n$ in the graph) from the initial point. This condition was not needed in
the much weaker result of Frieze, Kleinberg, Ravi and Debany~\cite{FKRD},
but this is not surprising. The very simple analysis in~\cite{FKRD}
ignored these long-range steps of the branching random walk (as well
as other complications), which is why it
does not give the correct limiting size of the giant component, or indeed
the asymptotic critical probability.

\subsection{The case $n=\Theta(\omega)$}\label{noc}

In this subsection we briefly consider the emergence of the giant
component in $\Gnop$ or $\Gnop'$ when $\omega\to\infty$
and $n=\Theta(\omega)$. We assume throughout that $n\sim C\omega$
for some constant $C>0$. As before, the relevant normalization
is to take $p=p(n)=\la/\omega$, with $\la$ constant.

As noted above, Theorem~\ref{th_2} does not apply in this setting.
However, our methods still allow us to obtain the critical
point, and indeed the size of the giant component, in terms
of the solutions to certain non-linear equations. As this problem
is of rather specialized interest, we shall only outline the results and
arguments.

Let $\PP_\la$ be a Poisson process on $\RR$ with intensity $\la$.
For $d>1$, let $g_\la(d)$ be the probability that within the interval
$[0,d]$ there is an interval of length $1$ containing no points
of $\PP$, and let $r_\la(d)=1-g_\la(d)$. For $0\le d\le 1$ we set $g_\la(d)=0$ and $r_\la(d)=1$.
Given that both $0$ and $d$ are points of $\PP_\la$, then $r_\la(d)$ is
the probability that they are joined in the graph with vertex set
$\PP_\la$ and edges 
between all points at distance at most $1$. The function
$r_\la(d)$ will be key to our analysis.
(This function $r_\la(d)$ has been studied in other contexts. It is 
the probability that the largest gap between $\Po(\la d)$ points
uniformly and independently distributed in $[0,1]$ is at most $1/d$;
this largest gap is also a version of a scan statistic;
see Glaz, Naus and Wallenstein~\cite{ScanStatistics}, for example.
An explicit, but rather complicated, expression
for $r_\la(d)$ is easily derived from the formula given
by Stevens \cite{Stevens}
for the case of the largest gap between a fixed number of random
points on a circle.)

Let us start with the simpler toroidal case. Let $\PP_{\la,C}$ be a Poisson
process of intensity $\la$ on the circle $\TT_C=\RR/(C\Z)$ with circumference
$C$. Adding in the extra point $0$, let $Z=Z(\la,C)$ denote (the distribution of)
the number of other points of $\PP_{\la,C}$ that may be reached from $0$
in steps of size at most $1$. Let $\bp_{\la,C}$ be the branching process
in which the number of children of initial particle has the distribution
of the sum of two independent copies of $Z$, and the number of children
of each later particle has the distribution of $Z$.
Approximating the component exploration in $\Gnop'$ by the corresponding
branching random walk on the torus, the proof of Theorem~\ref{th_2} adapts
easily to show that if $p=\la/\omega$ and $n\sim C\omega$,
then the asymptotic fraction of vertices of $\Gnop'$
that are in the giant component is simply the survival
probability $\rho(\la,C)$ of $\bp_{\la,C}$.
In particular, the critical value $\lac=\lac(C)$ of $\la$ may be found by
solving $\E(Z(\la,C))=1$.

If the point $x\in \TT_{C}$ is present, then it may be reached
from $0$ with probability $1-g_\la(x)g_\la(C-x)$, so $\lac$ is given by
the solution to
\[
 \int_{0}^{C} \bb{1-g_\la(x)g_\la(C-x)} \dd x =1.
\]

The case of the grid is more interesting, but also more complicated.
In this case we must consider the restriction $\rw_{\la,C}$ of our original
branching random walk $\rw_\la$ to a square $[0,C]^2$.
Previously, we used a similar restriction with $C$ large as a tool in our analysis.
Here we are forced to analyze the restriction itself, and have no control over $C$.

Fixing $\la$ and $C$,
for $0\le x,y\le C$, let $\rhoh(x,y)$, $\rhov(x,y)$ and $\rho(x,y)$
denote respectively the survival probabilities in the cases where we start with one particle
of type `h' at $(x,y)$, one particle of type `v' at $(x,y)$, and one
particle of each type at $(x,y)$.
By symmetry we have $\rhoh(x,y)=\rhov(y,x)$. Also,
\[
 \rho(x,y)=1-(1-\rhoh(x,y))(1-\rhov(x,y)) = \rhoh(x,y)+\rhov(x,y)-\rhoh(x,y)\rhov(x,y).
\]
Starting with two particles, one of each type, at a point $(x,y)$
chosen uniformly from $[0,C]^2$, the survival probability $\rho$ satisfies
\[
 \rho = \int_0^C\int_0^C \rho(x,y) = 2\int_0^C\int_0^C \rhoh(x,y) - \int_0^C\int_0^C \rhoh(x,y)\rhoh(y,x),
\]
so to understand $\rho$ it suffices to understand $\rhoh(x,y)$.

It is not hard to convince oneself that the proof of Theorem~\ref{th_2} can be
modified to show that if $p=\la/\omega$ and $n\sim C\omega$, then
the asymptotic fraction of vertices of $\Gnop$ in the largest component
is simply $\rho=\rho(\la,C)$. Filling in the details may well require
considerable work, however: for example, to make the sprinkling argument
work one needs to show that $\rho(\la,C)$ is continuous in $\la$.
Let us omit the details and turn instead to the study of $\rho(\la,C)$.

It is easy to see that one can express $\rhoh(x,y)$ as the maximum solution
to a certain non-linear integral equation; as the details are not very illuminating,
we ignore the size of the giant component and
consider only the question of when it emerges, i.e., when $\rho(\la,C)>0$.
At this point, it is tempting to claim that arguments of the type given
in~\cite{BJRkernels}, for example, 
show that $\rho(\la,C)>0$ if and only if the corresponding linearized operator
has norm strictly greater than $1$. However, there are some complications:
in~\cite{BJRkernels}, we assumed that the linear operator was symmetric and compact;
here it is neither!
These complications seem to be connected to the interchange of $x$ and
$y$ coordinates, or, putting it another way, to the one-dimensional
action of the operator on the two-dimensional space $[0,C]^2$.
It turns out that one can simplify things greatly by
considering two steps of the branching process at once.

Starting with a point at $(x,y)$ (of type `h', say, although it doesn't matter),
and taking two steps in our restricted branching random walk,
the expected number of points in $[x'+dx']\times[y'+dy']$ that we reach is
easily seen to be $\la^2 r_{\la}(|x-x'|)r_{\la}(|y-y'|)\dd x'\dd y'$.
Let $T_2$ be the corresponding integral operator on $L^2([0,C]^2)$,
so
\[
 (T_2(f))(x',y') = \int_0^C\int_0^C \la^2 r_\la(|x-x'|)r_\la(|y-y'|)f(x,y)\dd x\dd y.
\]
The operator $T_2$ is very well behaved: the kernel is clearly
symmetric (with respect to swapping $(x,y)$ and $(x',y')$)
and, since $0\le r\le 1$, the operator is compact.
Thus arguments of the type given~\cite{BJRkernels} do show that $\rho(\la,C)>0$
if and only if $\norm{T_2}>1$.

At this point, we can return to one dimension:
the operator $T_2$ has a unique (up to normalization)
eigenfunction $\phi$ with eigenvalue the norm of $T_2$.
Let $T$ be the operator on $L^2([0,C])$ defined by
\[
 (T(f))(x') = \int_0^C \la r_{\la}(|x-x'|)f(x) \dd x.
\]
From the form of $T_2$, it is not hard
to check that $\phi(x,y)=\psi(x)\psi(y)$,
where $\psi$ is the eigenfunction of $T$ with
maximum eigenvalue. In particular, $\norm{T_2}=\norm{T}^2$.
Although we only sketched the details, one has $\rho(\la,C)>0$ if and only if
$\norm{T_2}>1$, so in this setting the critical value of $\la$
is given by the solution to $\norm{T}=1$, noting that $\norm{T}$ depends
both on $C$ and on $\la$.

Unfortunately it seems unlikely that one can find $\norm{T}$ explicitly;
this is to be expected. The conclusion is that we find the same
sort of connection between the giant component in a certain random graph
and the norm of an integral operator on $L^2([0,C])$ as seen in~\cite{BJRkernels}.
This is especially interesting as one cannot directly apply the results of~\cite{BJRkernels}.
Indeed, here, the approximating branching process/random walk
is {\em not} of the multi-type Poisson form considered there.

\section{Discussion and variants}

The line-of-sight model $\Gop$ is unusual in various ways.
There are several random graph models in which one can determine the
exact threshold for the emergence of a giant component in terms of the survival
probability of a suitable
a branching process -- the simplest example is $G(n,p)$; a very
general inhomogeneous model with this property is that
of~\cite{BJRkernels}.  In lattice percolation models, on the other
hand, as a rule one cannot determine the exact critical
probability. Penrose~\cite{Penrose} defined a `spread-out'
percolation model as follows: let $S$ be a fixed bounded symmetric set
in $\RR^2$. Given parameters $\omega$ and $p$, for every pair of
points $v$, $w$ of $\Z^2$, join them by an edge with probability $p$
if $v-w\in \omega S$, independently of all other such
pairs. (Penrose's model is somewhat more general, but this is the
essence.) For fixed $\omega$, this is a lattice percolation model, and
the critical probability cannot be found exactly.  However, Penrose
showed that as $\omega\to\infty$, the asymptotic form of the critical
probability can be found -- the `critical expected degree' tends to
$1$. In~\cite{BJRspread}, it was shown that this result follows easily
from those of~\cite{BJRkernels}. Given the homogeneous nature
of the model, the critical expected degree tending to $1$ shows
that asymptotically, cycles do not matter; this is far from the
case in the line-of-sight model.

In terms of the critical probability, the line-of-sight model $\Gop$
behaves differently from Penrose's model: even when $\omega$ is large,
$\Gop$ contains many short cycles, so the critical expected degree,
$4\log(3/2)$, is larger than $1$. There are two important
differences between the models. Firstly, $\Gop$
is a site percolation model (vertices are selected at random),
and Penrose's model is a bond percolation one (edges are selected at random).
Secondly, in $\Gop$ the geometric condition for joining vertices
involves scaling a set $S\subset \RR^2$ of measure zero. We shall examine the effect
of these differences separately, by considering two variants of $\Gop$.

\medskip
There is a natural bond percolation variant of $\Gop$.
Given $\omega$ and $p$, let $\bGop$ be the graph on $\Z^2$ defined
as follows. For every pair of vertices $v$, $w$ 
such that the line segment $vw$ has length at most $\omega$ and is
horizontal or vertical, join $v$ and $w$ with probability $p$, independently
of all other such pairs. We consider the limit
$\omega\to\infty$ with $p=\la/\omega$, where $\la$ is constant.
In this model the density of short cycles is very low. Thus, as
we explore the component of a vertex in the usual way, the number
of vertices reached at each step may be approximated by a Galton--Watson
branching process in which each particle has a Poisson number of children
with mean $4\la$.
Moreover, the location of these vertices may be approximated
by a corresponding branching random walk $\bbp_{4\la}$. The arguments above
for $\Gop$ carry over to this setting, showing that the critical
probability $\pc(\omega)$ is asymptotically $1/(4\omega)$,
and that if $p=\la/\omega$ and $\omega\to\infty$ then
the infinite component contains a fraction $\bphi(4\la)+o(1)$
of the vertices, where $\bphi(\mu)$ is the survival probability
of the branching process underlying $\bbp_\mu$, i.e., the maximal
solution to $\bphi(\mu)=1-e^{-\mu\bphi(\mu)}$.
One can also obtain corresponding results for finite graphs, showing that
when $\omega\to\infty$, they behave much like the
Erd\H os-R\'enyi (or Gilbert) model $G(n,4\la/n)$.

\medskip
Turning to our second variant, let $C_\eps$ be the cross with (four) arms
of length $1$ and width $2\eps$, so
\[
 C_\eps = [-1,1] \times [-\eps,\eps] \ \cup\ [-\eps,\eps]\times [-1,1].
\]
Let $\Gope$ be the random graph defined as follows: define a graph
$G$ on $\Z^2$ by joining two vertices $v$ and $w$ if $v-w\in \omega C_\eps$.
Then select vertices independently with probability $p$ to form $\Gope$.
When $\eps=0$, this graph is exactly $\Gop$. For $\eps>0$, the area
of $\omega C_\eps$ is proportional to $\omega^2$, so we consider the limit
with $p=\la/\omega^2$ and $\omega\to\infty$.
In this context it is natural to rescale the vertex set, selecting
vertices of $(1/\omega)\Z^2$ with probability $\la/\omega^2$
and joining them if their difference lies in $C_\eps$.
In the limit, the selected vertices form a Poisson process on $\RR^2$
with intensity $\la$. This gives us another natural random graph
model, $\Gle$. The vertex set of $\Gle$ 
is a Poisson process $\PP_\la$ on $\RR^2$ with intensity $\la$,
and two vertices $v$ and $w$ are joined if and only if $v-w\in C_\eps$.
It is not hard to see that as $\omega\to\infty$ with $\eps$ fixed,
$\omega^2$ times the critical probability for percolation in $\Gope$
approaches the critical density $\lac$ for percolation
in $\Gle$.

The random graph $\Gle$ we have just defined is a special case
of a percolation model introduced by Gilbert~\cite{G61} in 1961,
Gilbert's {\em disc model}. To define this, let $S$
be a symmetric (in the sense $S=-S$) set in $\RR^2$
with Lebesgue measure $0<\area(S)<\infty$, and let
$G(\la,S)$ be the random graph whose vertex set is the Poisson
process $\PP_\la$, in which vertices are joined if their
difference lies in $S$, so $\Gle=G(\la,C_\eps)$.
Gilbert asked the following question: given $S$, what is the critical
value $\lac=\lac(S)$ above which $G(\lac,S)$ contains
an infinite component? Since $G(\la,S)$ and $G(\la/c^2,cS)$
have the same distribution as graphs, one often works
with the `critical expected degree', or {\em critical area}
$\ac(S)=\lac(S)\area(S)$, instead. (The term critical area
is used because of the standard normalization $\la=1$.)

As with many such percolation questions, it seems
impossible to determine $\ac(S)$ exactly.
The most studied cases are when $S$ is a disc $D$,
or when $S$ is a square $[-1,1]^2$.
The best rigorous bounds for $\ac(D)$, the
bounds $2.184\le \ac(D)\le 10.588$ proved by Hall~\cite{Hall85},
are not much better than Gilbert's original bounds; although
these bounds are
far apart, it seems hard to improve them significantly.
In practice, however, $\ac(D)$ is known quite precisely;
Balister, Bollob\'as and Walters~\cite{BBW} {\em proved}
that $[4.508,4.515]$ is a $99.99\%$ confidence interval for $\ac(D)$,
and $[4.392,4.398]$ for $\ac([-1,1]^2)$. See~\cite[Section 8.1]{BRbook}
for more details, as well as many references to heuristic bounds.

Recall that we introduced $\Gle=G(\la,C_\eps)$ as a variant of $\Gop$. 
In this context it is natural to take the limit $\eps\to 0$,
and ask whether the critical average degree
converges to the corresponding value for $\Gop$.
Since the area of the cross $C_\eps$ is $8\eps-4\eps^2$,
one might expect that
\begin{equation}\label{guess}
 8\eps\lac(C_\eps) \sim 4\log(3/2)
\end{equation}
as $\eps\to 0$. We shall see that this is not the case.
In fact, rather surprisingly, we can describe the
limit of $\eps\lac(C_\eps)$ in terms of the Gilbert model
with $S$ a square!

Given $\la>0$, let $G$ be the graph $G(\la,[-1,1]^2)$ conditioned
on the origin being a vertex. In other words, $G$ is the graph
on $\PP_\la\cup\{(0,0)\}$ in which we join two vertices
if their $\ell_\infty$ distance
is at most $1$.
Let
\[
 f_2(\la) = \E(|C_0|-1),
\]
where $C_0$ is the component of $G$ containing $(0,0)$, and
$|C_0|$ is the number of vertices in $C_0$. Note that we subtract
$1$ to avoid counting the `extra' vertex at the origin. Later
we shall consider an analogous quantity defined in terms
of a $d$-dimensional form of Gilbert's model; this is
the reason for the notation $f_2$.

\begin{theorem}\label{th_3}
The critical densities $\lac(C_\eps)$ satisfy
$
 \eps\lac(C_\eps) \to \mu
$
as $\eps\to 0$, where $\mu$ is the unique solution to
$
 f_2(\mu) = 1.
$
\end{theorem}
\begin{proof}
Note first that the definition of $\mu$ makes sense:
we have $f_2(0)=0$, while $f_2(\la)$ diverges at
$\lac([-1,1]^2)<\infty$.
Furthermore, up to this point $f_2(\la)$ is increasing and continuous;
hence there is a unique solution to $f_2(\mu)=1$.

Fix $\la>0$ and $\eps>0$; we shall take $\eps^{-1}\la$ as our density
parameter.

Let us condition on $v_0=(0,0)$ being a vertex of $\Glee=\GGlee$,
and explore the component of $v_0$ in this graph. As in our study
of $\Gop$, we break the exploration down into steps; in each step,
we explore either horizontally or vertically from a vertex $v$.
As before, we set $Y_0=\{v_0\}$, explore horizontally and vertically
from $v_0$, writing $Y_1$ for the set of vertices reached,
and then explore horizontally from those vertices of $Y_t$ reached during
vertical explorations, and {\em vice versa}. We write $Y_{t+1}$
for the new vertices found during such explorations,
so the vertex set of the component containing $v_0$ is the disjoint union of the sets
$Y_t$, $t\ge 0$. This time, however,
the horizontal and vertical explorations are a little more complicated.

Let us say that an edge $vw$ of $\Glee$ is {\em horizontal} if $v-w\in
[-1,1]\times [-\eps,\eps]$.  The definition of {\em vertical} edges is
analogous. Note that some edges are both horizontal and vertical, but
the proportion of such edges tends to $0$ as $\eps\to 0$.

To explore horizontally from $v$, let $H_v$
consist of all new (not reached in previous explorations)
vertices $v'$ of $\Glee$ joined to $v$ by horizontal edges,
together with all new vertices $v''$ joined to such vertices $v'$
by horizontal edges, and so on.
In other words, writing $\Glee^\mathrm{h}$ for the subgraph of $\Glee$
formed by all horizontal edges, $H_v\cup\{v\}$ is simply
the component containing $v$ in the subgraph of $\Glee^\mathrm{h}$
induced by the $v$ and the new vertices.
In particular, if this is our first exploration,
then $v=v_0=(0,0)$, and $H_v\cup\{v\}$ is simply the component
$C_0$ of the origin in $\Glee^\mathrm{h}$.

Now $\Glee^\mathrm{h}$ is exactly the Gilbert graph
$G(\la/\eps,[-1,1]\times[-\eps,\eps])$, conditioned on $(0,0)$ being a vertex.
Scaling vertically by a factor $\eps$, the distribution of this
random graph (as an abstract graph) is identical to that of 
$G(\la,[-1,1]^2)$.
Let $\DD_2(\la)$ denote the distribution of $|C_0|-1$, where
$C_0$ is the component of the origin in this graph,
so $f_2(\la)=\E(\DD_2(\la))$ by definition.
Then $|H_{v_0}|$ has the distribution $\DD_2(\la)$.

In later steps of the exploration, the restriction to new vertices
ensures that the distribution of $H_v$, given the exploration
so far, is stochastically dominated by that of $H_{v_0}$,
i.e., by $\DD_2(\la)$. The same holds for vertical explorations.
It follows that the sequence $|Y_0|,|Y_1|,\ldots$ is stochastically
dominated by $|X_0|,|X_1|,\ldots$, where $(X_t)$ is a branching process
in which the number of children of each particle except the initial
one has the distribution $\DD_2(\la)$, and different particles
have children independently. The number of children of the initial
particle is distributed as the sum of two independent copies of $\DD_2(\la)$, 
since we explore horizontally and vertically from $v_0$.

Since $\E(\DD_2(\la))=f_2(\la)$, this shows that
$\Glee$ has no giant component when $\la<\mu$, where
$\mu$ satisfies $f_2(\mu)=1$. Hence, $\eps\lac(C_\eps)\le \mu$.

\medskip
It remains to show that if $\la>\mu$ is constant, then for
$\eps$ small enough we have percolation in the graph $\Glee$.
The argument is very similar to that for $\Gop$, so we give only
the briefest outline, emphasizing the differences.
In doing so we may assume that $\la< \lac([-1,1]^2)$:
if $\la>\lac([-1,1]^2)$ then the horizontal subgraph of
$\Glee$ has the same distribution as a graph as
the supercritical Gilbert graph $G(\la,[-1,1]^2)$,
so this subgraph already contains an infinite component.

When exploring horizontally in $\Glee$, there is {\em a priori}
no bound on how far we may drift vertically. To deal with this,
let us choose a large constant $A$ (depending only on $\la$),
and limit each horizontal or vertical exploration to $A$ steps.
Let $\DD_2^{(A)}(\la)$ be the distribution of the number of points
reached from $v_0$ in this restricted horizontal exploration.
Since the horizontal subgraph of $\Glee$ is subcritical,
the distribution $\DD_2^{(A)}(\la)$ converges to $\DD_2(\la)$
as $A\to\infty$. In particular, since $\la>\mu$, we may choose
$A$ so that $\E\bb{\DD_2^{(A)}(\la)}> 1$.

With this modified exploration, the coupling argument
used in $\Gop$ carries over to show that we may couple
the first $O(1)$ steps of the exploration with the appropriate
branching process with an error probability of $O(\eps)$.
The main difference is that rather than claiming lines,
new vertices we find claim horizontal or vertical strips
of width $2A\eps$. As before, we can in fact couple
with a supercritical branching random walk, and use oriented percolation
to find an infinite component.
\end{proof}

Let $A_\eps$ be the annulus centered on the origin with
outer radius $1$ and inner radius $1-\eps$.
Independently, Franceschetti, Booth, Cook,
Meester and Bruck~\cite{FBCMB},
and Balister, Bollob\'as and Walters~\cite{BBW_ann} showed
that $\lim_{\eps\to 0}\ac(A_\eps)=1$. This result is similar
in spirit to the results here: one shows local approximation by
a branching process (in this case, a very simple process,
since the graph contains few short cycles), and then must work
to deduce percolation in the supercritical case.

Let $S_\eps$ be the `square annulus' $[-1,1]^2\setminus [-(1-\eps),1-\eps]^2$.
One might expect $\lim_{\eps\to 0}\ac(S_\eps)$ to also equal $1$,
but this is not true: in~\cite{BBW_ann} the bound
$\ac(S_\eps)\ge 1.014$ is proved for every $\eps>0$.
Our methods here show that $\ac(S_\eps)$ does converge as $\eps\to 0$,
and give a description of the limiting value.

Let $F=[-1,1]^2\times \{-1,1\}$, so $F\subset \RR^2\times \Z$ 
consists of two opposite faces of the cube $[-1,1]^3$.
We define a variant $G(\la,F)$ of Gilbert's model as follows:
the vertex set is a Poisson process of intensity $\la$ on $\RR^2\times \Z$,
with the origin added, and two vertices $v$ and $w$ are joined
if and only if $v-w\in F$. Let $C_0$ be the component of
the origin in $G(\la,F)$, and let $f(\la)=\E(|C_0|-1)$.

\begin{theorem}\label{th_A}
The critical area $\ac(S_\eps)$ tends to $2\mu$ as $\eps\to 0$, 
where $\mu$ is the unique solution to $f(\mu/8)=1$.
\end{theorem}
Note that the expected degree in $G(\la,F)$ is $8\la$; thus the limiting
expected area or degree in \refT{th_A} is twice the expected degree at which
$\E(|C_0|-1)=1$ holds in $G(\la,F)$. In this sense, \refT{th_A} is analogous
to \refT{th_3}. 
\begin{proof}
The proof is very similar to that of \refT{th_3}, so we omit the details.
The key observation is as follows: call an edge $vw$ of $G(\la,S_\eps)$
`horizontal' if $v-w$ lies in one of the vertical sides of the annulus,
i.e., if $\pm(v-w)\in [1-\eps,1]\times [-1,1]$, and define
vertical edges similarly. As $\eps\to 0$,
it is easy to check that almost all short cycles
in $G(\la/\eps,S_\eps)$ consist entirely of horizontal edges, or entirely of
vertical edges. As before, we break down the neighbourhood exploration
process in $G(\la/\eps,S_\eps)$ into horizontal and vertical explorations.
In a horizontal exploration, say, we look for all new vertices
that can be reached by horizontal edges.
As before, the asymptotic condition for criticality
is that the expected
number of new vertices found in a single horizontal (or vertical)
exploration is $1$. Rescaling horizontally by $2/\eps$,
a horizontal exploration may be coupled with an exploration of
the component of the origin
in the Gilbert model $G(\la/2,T_\eps)$, where $T_\eps$ consists of two squares
of side $2$ centered at $(\pm(2\eps^{-1}-1),0)$.
If we take $k\le 1/\eps$ steps from the origin, with the displacement
of each step in $T_\eps$, then we can tell from our final position
how many steps were to the right and how many to the left. It follows
that for the first $1/\eps$ steps, the explorations in $G(\la/2,T_\eps)$
and in $G(\la/2,F)$ may be regarded as identical.
It follows that $\eps\lac(S_\eps)\sim 2(\mu/8)$, where $f(\mu/8)=1$.
Since the area of $S_\eps$ is asymptotically $8\eps$, the result follows.
\end{proof}

Since it is much easier to estimate (by simulation)
the point at which an expectation crosses $1$ than the point at which
it diverges, \refT{th_A} makes it much easier to estimate
$\lim_{\eps\to 0}\ac(S_\eps)$; a simple simulation suggests
that the limit is $1.11406\pm 0.00001$.

\bigskip
The models above make just as good sense in any dimension.
For $d\ge 1$ and $\la\ge 0$, let $G=G(\la,[-1,1]^d)$ be the graph whose
vertex set consists of a Poisson process in $\RR^d$ of intensity $\la$
together with the origin, in which two vertices are joined if they
are within $\ell_\infty$ distance $1$. Let
\[
 f_d(\la) = \E(|C_0|)-1,
\]
where $C_0$ is the component of $G$ containing the origin,
so $f_d(\la)=\infty$ if $\la\ge\lac$, where $\lac$ is the critical
density for $G(\la,[-1,1]^d)$. In general, one cannot hope
to evaluate $f_d(\la)$ exactly. However,
when $d=1$, $|C_0|-1$ has exactly the distribution $\Ga_\la^{(2)}$
defined earlier, so $f_1(\la)=2(e^{\la}-1)$,
and $f_1(\log(3/2))=1$. It is easy to check that $f_2(\la/2)>f_1(\la)$
for any $\la>0$. (This is the natural comparison,
as the expected degree of a vertex of the relevant graph is $2^d\la$.)
It follows that the quantity $\mu$ defined in \refT{th_3}
is strictly less than $\log(3/2)/2$, so \refT{th_3}
shows that our `guess' \eqref{guess} does not hold.
(A quick simulation suggests that $\mu$ is around 0.177635.) 
In other words, the critical expected degree in models defined using
a cross with very thin arms does not approach that of a model
defined using a cross whose arms have width $0$. (Compare
$8\times 0.177635=1.421\dots$ with $4\log(3/2)=1.621\dots$.)

\bigskip
We finish by turning to one final class of variants of $\Gop$, namely
the natural generalizations to higher dimensions. 
Of course, there is more than one natural equivalent of a $2$-dimensional cross
in higher dimensions: a union of line segments, or a union of $(d-1)$-dimensional
hypercubes. Fortunately one can cover both generalizations 
in a single definition.

Given $d\ge 2$, $1\le r\le d-1$, $\omega\ge 1$ and $0<p<1$, 
let $\Gdrop$ be the random graph defined as follows:
for the vertex set, select points of $\Z^d$ independently
with probability $p$. Join two vertices $v$ and $w$ if and only
if $v$ and $w$ differ in at most $r$ coordinates, and differ
in each of these coordinates by at most $\omega$.
If $d=2$ and $r=1$ then we obtain $\Gop$.

\begin{theorem}\label{th_4}
Let $\pc(d,r,\omega)$ denote the critical probability $p$
for percolation in $\Gdrop$. Then
\[
 \omega^r \pc(d,r,\omega) \to \la_{d,r}
\]
when $\omega\to\infty$ with $d\ge 2$
and $1\le r\le d-1$ fixed,
where $\la_{d,r}$ is the unique solution to
\[
 \left(\binom{d}{r}-1\right) f_r(\la_{d,r}) =1.
\]
\end{theorem}
Recalling that $f_1(\la)=2(e^{\la}-1)$, we see that
\[
 \la_{d,1} = \log\left(\frac{2d-1}{2d-2}\right),
\]
so \refT{th_4} generalizes the asymptotic critical probability result
for $\Gop$ given by \refT{th_1}. Note that the function $f_2(\mu)$,
which appeared in \refT{th_3} in the analysis of a certain $2$-dimensional graph,
appears here in the analysis of
the very different graphs $G_{d,2,\omega,p}$, $d\ge 3$.

The proof of \refT{th_4} is very similar to that of \refT{th_1}, so we omit
it. The main difference is that instead of exploring horizontally or vertically,
i.e., always exploring within a line, we always explore within an
affine 
subspace of dimension
$r$ generated by $r$ of the coordinate axes. The factor $\binom{d}{r}-1$
above appears because there are $\binom{d}{r}$ such subspaces
through any given point, and we reach any vertex
other than the initial one by exploring within one of them.
This time, when we find a vertex $v$ of $\Gdrop$
by exploring within a certain subspace, this vertex immediately claims
the other $\binom{d}{r}-1$ subspaces in which it lies. As before, we omit
testing points $w$ in two (or more) claimed subspaces. When exploring
from $v$ in a subspace $S$, any such point $w$ lies in a subspace $S'$
through a vertex $v'$ of $C_0$ that we have already found, with $S'$ and $S$
unequal and therefore not parallel. Since $S$ and $S'$ intersect
in a space of lower dimension, in the early stages of our exploration
almost all points of $S$ are claimed only by $v$, and the coupling
goes through as before.

Just as for $\Gop$, we could `thicken' the generalized crosses defining $\Gdrop$
by $\eps$ and let $\eps\to 0$, obtaining a generalization of \refT{th_3}.
In this case the limit $\mu$ in \refT{th_3} is replaced by the solution to
$
 \left(\binom{d}{r}-1\right) f_d(\mu) =1.
$

\bigskip\noindent
{\bf Acknowledgement}
This paper was inspired by Alan Frieze's talk on his joint work with
Kleinberg, Ravi and Debany at the Oberwolfach meeting `Combinatorics,
Probability and Computing', November 2006.

\end{document}